# DESIGN OF FULL ORDER PROPORTIONAL-INTEGRAL OBSERVER FOR DISCRETE-TIME LINEAR TIME-INVARIANT SYSTEMS

Konstadinos H. Kiritsis

Hellenic Air Force Academy, Department of Aeronautical Sciences, Division of Automatic Control, Dekelia Air Base, PC 13671 Acharnes, Attikis,Tatoi, Greece
e-mail: konstantinos.kyritsis@hafa.haf.gr

**Abstract**

*This paper is devoted to the design of full order proportional-integral observer for the state estimation of discrete-time linear time-invariant systems. In particular, explicit necessary and sufficient conditions are established for the existence of proportional-integral observer for the state estimation of discrete-time linear time-invariant systems and a simple procedure is given for the construction of the observer. Our approach is based on properties of real and polynomial matrices.*

## 1. INTRODUCTION

In 1971 Luenberger proposed the full order observer for the state estimation of linear time-invariant systems [1]. In [2] Wojciechowski added an additional term to Luenberger's full order observer for the state estimation of single-input single-output linear time-invariant systems. This term is proportional to the integral of the output estimation error. The resulting new observer was called proportional-integral observer and has a long and rich history. The main results of [2] were generalized to linear multivariable time-varying systems in [3], in particular in [3] implicit necessary and sufficient conditions for the existence of proportional-integral observer for the state estimation of linear multivariable time-varying systems have been established. In [4] a reduced order proportional-integral observer for the state estimation of linear multivariable time-varying systems was first considered. In [5] the robustness property of feedback control systems using a proportional-integral observer was studied. In [6] necessary and sufficient conditions have been derived under which the proportional-integral observer achieves Exact Loop Transfer Recovery for continuous-time linear time-invariant systems. Similar results have been obtained in [7] for discrete-time linear time-invariant systems. In [8] it was proved that the proportional-integral observer can estimate the state not only of linear time-invariant systems but also of systems with arbitrary external input which appear as unknown input, nonlinearity or unmodeled dynamics.

In [9] it was shown that, for some classes of systems the proportional-integral observer has the ability to completely decouple the modeling uncertainties while keeping satisfactory convergence properties. Furthermore a comparison of classical proportional observer to proportional-integral observer was given using a simulation example. A parametric eigenstructure assignment design approach for proportional-integral observers for the state estimation of continuous-and discrete-time linear time-invariant systems was proposed in [10] and [11] respectively. In [12] a proportional-integral observer based sliding mode controller was proposed for nonlinear hydraulic differential cylinder systems affected by uncertainties. In [13] an optimization method based on a genetic algorithm for the computation of gains of proportional-integral observer for the estimation of state variables of an induction motor is presented. Proportional-integral observer-based approaches for fault detection were developed in [14-16] and references given therein. The proportional-integral observer literature is extremely rich; for more complete references, we refer the reader to [17], [18] and [19]. To the best of our knowledge the problem of design of full order proportional-integral observer for the state estimation of discrete-time linear time-invariant systems, is still an open problem. This motivates the present study. Associated with the design of full order proportional-integral observer for the state estimation of discrete-time linear time-invariant systems are two fundamental questions, i.e. the question of solvability and the question of computability. A major effort in solvability is to determine necessary and sufficient conditions for the existence of a full order proportional-integral observer. The main concern associated with computability, on the other hand, is to develop a procedure for the construction of the proportional-integral observer.

In this paper, these questions have been completely answered. In particular, by using basic concepts and basic results from linear systems and control theory as well as of the theory of matrices are established explicit necessary and sufficient conditions for the existence of a full order proportional-integral observer for the state estimation of discrete-time linear time-invariant systems and a simple procedure is given for the construction of the proportional-integral observer.

## 2. PROBLEM STATEMENT

Consider a discrete-time linear time-invariant system described by the following state-space equations

$$\mathbf{x}(k+1) = \mathbf{A}\mathbf{x}(k) + \mathbf{B}\mathbf{u}(k) \qquad (1)$$

$$\mathbf{y}(k) = \mathbf{C}\mathbf{x}(k) \qquad (2)$$

where $\mathbf{A}$, $\mathbf{B}$ and $\mathbf{C}$ are real matrices of size ($n \times n$), ($n \times m$) and ($p \times n$) respectively, $\mathbf{x}(k)$ is the state vector of size ($n \times 1$), $\mathbf{u}(k)$ is the vector of inputs of size ($m \times 1$) and $\mathbf{y}(k)$ is the vector of outputs of size ($p \times 1$). In what follows without any loss of generality we assume that

$$rank[\mathbf{C}] = p \qquad (3)$$

Let us consider a discrete-time linear time-invariant system described by the equations

$$\hat{\mathbf{x}}(k+1) = (\mathbf{A} - \mathbf{LC})\,\hat{\mathbf{x}}(k) + \mathbf{L}\mathbf{y}(k) + \mathbf{B}\mathbf{u}(k) + \mathbf{F}\mathbf{v}(k) \qquad (4)$$

$$\mathbf{v}(k+1) = \mathbf{v}(k) + [\mathbf{y}(k) - \mathbf{C}\hat{\mathbf{x}}(k)] \qquad (5)$$

where $\hat{\mathbf{x}}(k)$ is the state vector of dimensions ($n \times 1$), $\mathbf{v}(k)$ is a vector of size ($p \times 1$) and $\mathbf{L}$ and $\mathbf{F}$ are real matrices of size ($n \times p$) respectively. The discrete-time linear time-invariant system described by the equations (4) and (5) is a proportional-integral observer of order $n$ for the system described by the equations (1) and (2), if and only if for arbitrary initial conditions $\hat{\mathbf{x}}(0)$, $\mathbf{x}(0)$ and any input $\mathbf{u}(k)$, the following relationships hold [11]

$$\lim_{k \to +\infty} \mathbf{e}(k) = 0 \qquad (6)$$

$$\lim_{k \to +\infty} \mathbf{v}(k) = 0 \qquad (7)$$

where $\mathbf{e}(k) = [\hat{\mathbf{x}}(k) - \mathbf{x}(k)]$ is the state estimation error, $\hat{\mathbf{x}}(k)$ is an estimate of the state vector $\mathbf{x}(k)$ and $\mathbf{v}(k)$ is a vector representing the integral of the weighted output estimation error [11]. The relationships (6) and (7) are simultaneously satisfied if and only if the matrix

$$\begin{bmatrix} \mathbf{A} - \mathbf{LC} & \mathbf{F} \\ -\mathbf{C} & \mathbf{I}_p \end{bmatrix} \qquad (8)$$

of size (($n+p$) x ($n+p$)) is Schur stable, i.e. all its eigenvalues have magnitude less than 1 [11]. Thus the problem of the design of the proportional-integral observer of order $n$ can be stated as follows: Do there exist real matrices $\mathbf{L}$ and $\mathbf{F}$ of appropriate dimensions such that the matrix given by (8) is Schur stable? If so, give conditions for existence and a procedure for the calculation of the real matrices **L** and **F**.

## 2. BASIC CONCEPTS AND PRELIMINARY RESULTS

This section contains lemmas, which are needed to prove the main results of this paper and some basic notions from linear systems and control theory as well as of the theory of matrices that are used throughout the paper. Let $\mathcal{R}$ be the field of real numbers. Also let $\mathcal{R}[z]$ be the ring of polynomials with coefficients in $\mathcal{R}$. Further, let $C$ be the field of complex numbers, also let $C^+$ be the set of all complex numbers $\lambda$ with $|\lambda| \geq 1$. The units of $\mathcal{R}[z]$ are polynomials of zero degree, i.e. all nonzero finite real numbers. A polynomial over $\mathcal{R}[z]$ is said to be non-unit if and only if it has nonzero degree. A matrix whose elements are polynomials over $\mathcal{R}[z]$ is termed a polynomial matrix. A polynomial matrix $\mathbf{U}(z)$ over $\mathcal{R}[z]$ of size *(q x q)* whose determinant is a unit of $\mathcal{R}[z]$ is termed unimodular matrix [20]. Every polynomial matrix $\mathbf{M}(z)$ of size *(m x p)* with $rank[\mathbf{M}(z)]=r$, can be expressed as [20]

$$\mathbf{U}_1(z)\,\mathbf{M}(z)\,\mathbf{U}_2(z) = \begin{bmatrix} \mathbf{M}_r(z) & \mathbf{0} \\ \mathbf{0} & \mathbf{0} \end{bmatrix} \quad (9)$$

The polynomial matrices $\mathbf{U}_1(z)$ and $\mathbf{U}_2(z)$ are unimodular and the non-singular polynomial matrix $\mathbf{M}_r(z)$ of size *(r x r)* in (9) is given by

$$\mathbf{M}_r(z) = \text{diag}\,[a_1(z),\,a_2(z),\,....,\,a_r(z)] \quad (10)$$

The nonzero polynomials $a_i(z)$ for $i=1,2,...,r$ are termed invariant polynomials of $\mathbf{M}(z)$ and have the following property

$$a_i(z) \text{ divides } a_{i+1}(z), \text{ for } i=1,2,...,r\text{-}1 \quad (11)$$

The relationship (9) with $\mathbf{M}_r(z)$ given by (10) is called Smith-McMillan form of $\mathbf{M}(z)$ over $\mathcal{R}[z]$. Since the matrices $\mathbf{U}_1(z)$ and $\mathbf{U}_2(z)$ are unimodular and the polynomial matrix $\mathbf{M}_r(z)$ given by (10) is non-singular, from (9) and (10) it follows that

$$rank[\mathbf{M}(z)] = rank[\mathbf{M}_r(z)] = r \quad (12)$$

**Definition 1:** The nonzero polynomial $c(z)$ over $\mathcal{R}[z]$ is said to be strictly Schur if and only if $c(z) \neq 0,\ \forall z \in C^+$.

**Definition 2:** Matrix **A** over $\mathcal{R}$ of size *(n x n)* is said to be Schur stable if and only if all eigenvalues of the matrix **A** have magnitude less than 1, or alternatively if and only if the characteristic polynomial of the matrix **A** is a strictly Schur polynomial

**Definition 3:** Let **A** and **C** be matrices over $\mathcal{R}$ of size *(n x n)* and *(p x n)*, respectively. Then the pair (**A**, **C**) is said to be detectable if and only if there exists a matrix **K** over $\mathcal{R}$ of size *(n x p)* such that the matrix [**A**+**KC**] is Schur stable [21].

**Definition 4**: Let **A** and **C** be matrices over $\mathcal{R}$ of size *(n x n)* and *(p x n)*, respectively and **C** not zero. Then an eigenvalue $\lambda$ of the matrix **A** is said to be observable, if and only if the following condition holds [22]:

$$rank \begin{bmatrix} \mathbf{C} \\ \mathbf{I}_n \lambda - \mathbf{A} \end{bmatrix} = n$$

Let **A** be a real matrix of size *(n x n)*. The spectrum of the matrix **A**, is the set of all its eigenvalues and is denoted by $\sigma(\mathbf{A})$. An eigenvalue $\lambda$ of **A** is called a stable eigenvalue if and only if $|\lambda| < 1$. The eigenvalue $\lambda$ of the matrix **A** is said to be unstable if and only if $|\lambda| \geq 1$.

**Lemma 1:** Let **A** and **C** be matrices over $\mathcal{R}$ of size *(n x n)* and *(p x n)*, respectively and **C** not zero. Further let $\sigma(\mathbf{A})$ be the spectrum of the matrix **A**. The pair (**A**, **C**) is detectable if and only if one of the following equivalent conditions holds [23]:

(a) $rank \begin{bmatrix} \mathbf{C} \\ \mathbf{I}_n z - \mathbf{A} \end{bmatrix} = n$, $\forall z \in C^+$

(b) $rank \begin{bmatrix} \mathbf{C} \\ \mathbf{I}_n \lambda - \mathbf{A} \end{bmatrix} = n$, $\forall \lambda \in \sigma(\mathbf{A})$ with $|\lambda| \geq 1$

From condition (b) of Lemma 1 it follows that the pair (**A**, **C**) is detectable if and only if all unstable eigenvalues of the matrix **A** are observable [23].

**Lemma 2.** Let **A** and **C** be matrices over $\mathcal{R}$ of size *(n x n)* and *(p x n)*, respectively. Then the pair (**A**, **C**) is observable if and only if for every monic polynomial $c(z)$ over $\mathcal{R}[z]$ of degree *n* there exists a matrix **K** over $\mathcal{R}$ of size *( n x p)*, such that the matrix [**A**+**KC**] has characteristic polynomial $c(z)$[20].

The standard decomposition of unobservable systems given in the following Lemma was first published by Kalman in [24] and can be also found in any standard text of linear systems theory.

**Lemma 3:** Let **A** and **C** be matrices over $\mathcal{R}$ of size *(n x n)* and *(p x n)*, respectively and **C** not zero. Further, let the pair (**A**, **C**) is unobservable. Then there exists a non-singular matrix **T** of size *(n x n)* such that

$$\mathbf{T}^{-1}\mathbf{AT} = \begin{bmatrix} \mathbf{A}_{11} & \mathbf{0} \\ \mathbf{A}_{21} & \mathbf{A}_{22} \end{bmatrix}$$

$$CT = [C_1, 0]$$

The pair $(A_{11}, C_1)$ is observable and the eigenvalues of the matrix $A_{22}$ are the unobservable eigenvalues of the pair $(A, C)$.

**Lemma 4:** Let $A$ and $C$ be matrices over $\mathcal{R}$ of size $(n \times n)$ and $(p \times n)$, respectively and $C$ not zero. Further let

$$A = T \begin{bmatrix} A_{11} & 0 \\ A_{21} & A_{22} \end{bmatrix} T^{-1}, \quad C = [C_1, 0] T^{-1}$$

with $(A_{11}, C_1)$ observable. If the pair $(A, C)$ is detectable then the matrix $A_{22}$ is Schur stable [23].

**Lemma 5:** Let $A$ be a matrix over $\mathcal{R}$ of size $(n \times n)$. Then the matrix $A$ is Schur stable if and only if the following condition holds:

(a) $rank[I_n z - A] = n$, $\forall z \in C^+$

*Proof:* Let $A$ be a strictly Schur matrix over $\mathcal{R}$. From Definition 2 it follows that the characteristic polynomial $c(z)$ of the matrix $A$ is a strictly Schur polynomial and therefore from Definition 1 it follows that

$$c(z) \neq 0, \forall z \in C^+ \tag{13}$$

The Smith-McMillan form of polynomial matrix $[I_n z - A]$ over $\mathcal{R}[z]$ is given by

$$K(z) [I_n z - A] L(z) = [diag[c_1(z), c_2(z), ...., c_n(z)]] \tag{14}$$

where $K(z)$ and $L(z)$ are unimodular matrices over $\mathcal{R}[z]$. The polynomials $c_i(z)$ for $i=1, 2,...,n$ are the invariant polynomials of the matrix $A$ and therefore their product is the characteristic polynomial $c(z)$ of the matrix $A$ [20], that is

$$c(z) = \prod_{i=1}^{n} c_i(z) \tag{15}$$

From (13) and (15) it follows that

$$c_i(z) \neq 0, \forall z \in C^+, \forall i=1,2,...,n \tag{16}$$

from (16) it follows that

$$rank\{diag[c_1(z), c_2(z), ...., c_n(z)]\} = n, \forall z \in C^+ \tag{17}$$

Since $K(z)$ and $L(z)$ are unimodular matrices over $\mathcal{R}[z]$, from (12) and (14) we obtain:

$$rank[I_n z - A] = rank\{diag[c_1(z), c_2(z), ...., c_n(z)]\} \tag{18}$$

Relationships (17) and (18) imply that

$$rank[I_n z - A] = n, \forall z \in C^+ \tag{19}$$

This is condition (a) of the Lemma. To prove sufficiency, we assume that condition (a) holds. Since by assumption condition (a) holds we have that

$$rank[\mathbf{I}_n z - \mathbf{A}] = n , \forall z \in C^+ \qquad (20)$$

From (18) and (20) we obtain:

$$rank\{diag[c_1(z), c_2(z), ...., c_n(z)]\} = n, \forall z \in C^+ \qquad (21)$$

From (21) it follows that

$$c_i(z) \neq 0 , \forall z \in C^+ , \forall i=1,2,...,n \qquad (22)$$

From (22) it follows that

$$\Pi_{i=1}^{n} c_i(z) \neq 0 , \forall z \in C^+ , \forall i=1,2,...,n \qquad (23)$$

Relationships (23) and (15) imply

$$c(z) \neq 0 , \forall z \in C^+ \qquad (24)$$

Relationship (24) and Definition 1 imply that $c(z)$ is a strictly Schur polynomial. Since by assumption $c(z)$ is the characteristic polynomial of the real matrix **A**, from Definition 2 it follows that the matrix **A** is Schur stable. This completes the proof.

The following Lemma is based on the results of [20].

**Lemma 6:** Let **A** and **C** be matrices over $\mathcal{R}$ of size *(n x n)* and *(p x n)*, respectively and **C** not zero. Further, let the pair (**A**, **C**) be detectable. Then there exists a matrix **K** over $\mathcal{R}$ of size *(n x p)* such that the matrix [**A**+**KC**] is Schur stable.

*Proof*: Let the pair (**A**, **C**) be detectable. Delectability of the pair (**A**, **C**) implies that the pair (**A**, **C**) is either observable or unobservable with stable unobservable eigenvalues. If the pair (**A**, **C**) is observable, then from Lemma 2 it follows that there exists a matrix **K** over $\mathcal{R}$ of appropriate dimensions such that

$$det[\mathbf{I}_n z - \mathbf{A} - \mathbf{KC}] = det[\mathbf{I}_n z - \mathbf{A} - \mathbf{KC}] = c(z) \qquad (25)$$

where $c(z)$ be an arbitrary monic, strictly Schur polynomial over $\mathcal{R}[z]$ of degree *n*. Since the notion of observability is a dual of reachability (i.e. observability of the pair (**A**, **C**) implies reachability of the pair(**A**$^T$, **C**$^T$)) [20], the matrix **K** can be calculated using known methods for the solution of pole assignment problem by state feedback [20]. Since $c(z)$ is the characteristic polynomial of the matrix [**A**+**KC**], from Definition 2 and (25) it follows that the matrix [**A**+**KC**] is Schur stable. If the pair (**A**, **C**) is unobservable with stable unobservable eigenvalues, then from Lemma 3 and Lemma 4 it follows that there exists a non-singular matrix **T** such that

$$\mathbf{A} = \mathbf{T}\begin{bmatrix}\mathbf{A}_{11} & \mathbf{0} \\ \mathbf{A}_{21} & \mathbf{A}_{22}\end{bmatrix}\mathbf{T}^{-1}, \quad \mathbf{C} = [\mathbf{C}_1, \mathbf{0}]\mathbf{T}^{-1} \qquad (26)$$

The pair ($A_{11}$, $C_1$) is observable and the matrix $A_{22}$ is Schur stable. Schur stability of the matrix $A_{22}$ and Definition 2 imply that the polynomial $\chi(z)$ given by

$$det[\ Iz - A_{22}] = \chi(z) \tag{27}$$

is a strictly Schur polynomial. Observability of the pair ($A_{11}$, $C_1$) and Lemma 2 imply the existence of a matrix $K_1$ over $\mathcal{R}$ of appropriate dimensions such that

$$det[\ Iz - A_{11} - K_1 C_1] = \varphi(z) \tag{28}$$

where $\varphi(z)$ is an arbitrary monic, strictly Schur polynomial over $\mathcal{R}[z]$ of appropriate degree. Since the notion of observability is a dual of reachability (i.e. observability of the pair ($A_{11}$, $C_1$) implies reachability of the pair($A_{11}^T, C_1^T$) ), the matrix $K_1$ can be calculated using known methods for the solution of pole assignment problem by state feedback [20]. Let

$$K = T \begin{bmatrix} K_1 \\ 0 \end{bmatrix} \tag{29}$$

Using (26) and (29) we have that

$$[A + KC] = = T \begin{bmatrix} A_{11} + K_1 C_1 & 0 \\ A_{21} & A_{22} \end{bmatrix} T^{-1} \tag{30}$$

while from (27), (28) and (30) we have that

$$det[(I_n z - A - KC] = \varphi(z)\chi(z) \tag{31}$$

Since by (27) and (28) the polynomials $\chi(z)$ and $\varphi(z)$ are strictly Schur, the polynomial $[\varphi(z)\chi(z)]$ is also a strictly Schur polynomial. Since by (31) the polynomial $[\varphi(z)\chi(z)]$ is the characteristic polynomial of the matrix $[A + KC]$, from Definition 2 it follows that the matrix $[A + KC]$ is Schur stable. This completes the proof.

## 4. MAIN RESULTS

The theorem that follows is the main result of this paper and gives the necessary and sufficient conditions for the existence of a full order proportional-integral observer for the state estimation of discrete-time linear time-invariant systems.

**Theorem 1.** The system described by equations (4) and (5) is a proportional-integral observer of order *n* of the system described by equations (1) and (2), if and only if the following condition holds:

(a) The pair (**A**, **C**) is detectable.

*Proof*: Let the system described by equations (4) and (5) is a proportional-integral observer of order *n* of the system described by equations (1) and (2). Then the real matrix of size $((n+p) \times (n+p))$ given by (8) is Schur stable. Schur stability of the matrix given by (8) and Lemma 5 imply

$$rank\begin{bmatrix} \mathbf{I}_n z - \mathbf{A} + \mathbf{LC} & -\mathbf{F} \\ \mathbf{C} & (z-1)\mathbf{I}_p \end{bmatrix} = (n+p) \ , \forall z \in C^+ \qquad (32)$$

Since the $(n+p)$ columns of the matrix on the left side of (32) are linearly independent over $C$, $\forall z \in C^+$, a subset of these columns consisting of the first $n$ columns must be also linearly independent over $C$, $\forall z \in C^+$; therefore

$$rank\begin{bmatrix} \mathbf{I}_n z - \mathbf{A} + \mathbf{LC} \\ \mathbf{C} \end{bmatrix} = rank\begin{bmatrix} \mathbf{C} \\ \mathbf{I}_n z - \mathbf{A} + \mathbf{LC} \end{bmatrix} = n \ , \forall z \in C^+ \qquad (33)$$

From (33) and after simple algebraic manipulations we obtain:

$$rank\begin{bmatrix} \mathbf{C} \\ \mathbf{I}_n z - \mathbf{A} + \mathbf{LC} \end{bmatrix} = rank\{\begin{bmatrix} \mathbf{I}_p & 0 \\ \mathbf{L} & \mathbf{I}_n \end{bmatrix}\begin{bmatrix} \mathbf{C} \\ \mathbf{I}_n z - \mathbf{A} \end{bmatrix}\} = n \ , \forall z \in C^+ \qquad (34)$$

Since matrix

$$\begin{bmatrix} \mathbf{I}_p & 0 \\ \mathbf{L} & \mathbf{I}_n \end{bmatrix} \qquad (35)$$

is non-singular, from (34) it follows that

$$rank\begin{bmatrix} \mathbf{C} \\ \mathbf{I}_n z - \mathbf{A} \end{bmatrix} = n \ , \forall z \in C^+ \qquad (36)$$

Relationship (36) and condition (a) of Lemma 1 imply that the pair $(\mathbf{A}, \mathbf{C})$ is detectable. This is condition (a) of the Theorem.

To prove sufficiency, we assume that condition (a) holds. Delectability of the pair $(\mathbf{A}, \mathbf{C})$ and Lemma 6 imply the existence a matrix $\mathbf{K}$ over $\mathcal{R}$ of size $(n \times p)$ such that the matrix $[\mathbf{A}+\mathbf{KC}]$ is Schur stable, that is

$$det[\mathbf{I}_n z - \mathbf{A} - \mathbf{KC}] = c(z) \qquad (37)$$

where $c(z)$ is a strictly Schur polynomial over $\mathcal{R}[z]$ of degree $n$. Matrix $\mathbf{K}$ in (37) can be calculated as in the proof of Lemma 6.

From (3) it follows that there exists a non-singular matrix $\mathbf{T}$ of size $(n \times n)$ such that

$$\mathbf{C} = [\mathbf{I}_p \ , 0]\mathbf{T} \qquad (38)$$

Let $\mathbf{\Phi}$ be an arbitrary nonzero Schur stable matrix over $\mathcal{R}$ of size $(p \times p)$. Furthermore, let $\mathbf{X}$ be a matrix over $\mathcal{R}$ of size $(n \times p)$ given by

$$\mathbf{X} = \mathbf{T}^{-1}\begin{bmatrix} (-\mathbf{\Phi} + \mathbf{I}_p) \\ \mathbf{\Lambda} \end{bmatrix} \qquad (39)$$

where $\mathbf{\Lambda}$ is an arbitrary matrix over $\mathcal{R}$ of size $((n-p) \times p)$. From (38) and (39) we have:

$$-\mathbf{CX} + \mathbf{I}_p = [-\mathbf{I}_p, 0]\mathbf{TT}^{-1}\begin{bmatrix} (-\mathbf{\Phi} + \mathbf{I}_p) \\ \mathbf{\Lambda} \end{bmatrix} + \mathbf{I}_p = \mathbf{\Phi} \qquad (40)$$

Now we form the matrix $\mathbf{M}$ over $\mathcal{R}$ of size $((n+p) \times (n+p))$ [25]

$$\mathbf{M} = \begin{bmatrix} \mathbf{I}_n & \mathbf{X} \\ \mathbf{0} & \mathbf{I}_p \end{bmatrix}$$

Matrix $\mathbf{M}$ over $\mathcal{R}$ is non-singular and its inverse is given by

$$\mathbf{M}^{-1} = \begin{bmatrix} \mathbf{I}_n & -\mathbf{X} \\ \mathbf{0} & \mathbf{I}_p \end{bmatrix}$$

We obtain:

$$\mathbf{M}^{-1} \begin{bmatrix} \mathbf{A} - \mathbf{LC} & \mathbf{F} \\ -\mathbf{C} & \mathbf{I}_p \end{bmatrix} \mathbf{M} = \begin{bmatrix} \mathbf{A} + (\mathbf{X} - \mathbf{L})\mathbf{C} & (\mathbf{A} - \mathbf{LC})\mathbf{X} + \mathbf{F} - \mathbf{X}(-\mathbf{CX} + \mathbf{I}_p) \\ -\mathbf{C} & -\mathbf{CX} + \mathbf{I}_p \end{bmatrix} \quad (41)$$

Furthermore, we set:

$$\mathbf{L} = \mathbf{X} - \mathbf{K} \quad (42)$$

$$\mathbf{F} = -(\mathbf{A} - \mathbf{LC})\mathbf{X} + \mathbf{X}(-\mathbf{CX} + \mathbf{I}_p) \quad (43)$$

Now by substituting (40), (42) and (43) into (41) we have:

$$\mathbf{M}^{-1} \begin{bmatrix} \mathbf{A} - \mathbf{LC} & \mathbf{F} \\ -\mathbf{C} & \mathbf{I}_p \end{bmatrix} \mathbf{M} = \begin{bmatrix} \mathbf{A} + \mathbf{KC} & \mathbf{0} \\ -\mathbf{C} & \mathbf{\Phi} \end{bmatrix} \quad (44)$$

Since by (37) matrix $[\mathbf{A}+\mathbf{KC}]$ is Schur stable and the nonzero real matrix $\mathbf{\Phi}$ is by assumption Schur stable, from Lemma 5 it follows that

$$rank[\mathbf{I}_n z - \mathbf{A} - \mathbf{KC}] = n, \forall z \in C^+ \quad (45)$$

$$rank[\mathbf{I}_p z - \mathbf{\Phi}] = p, \forall z \in C^+ \quad (46)$$

Relationships (45) and (46) imply

$$rank \begin{bmatrix} \mathbf{I}_n z - \mathbf{A} - \mathbf{KC} & \mathbf{0} \\ \mathbf{C} & \mathbf{I}_p z - \mathbf{\Phi} \end{bmatrix} = (n + p), \forall z \in C^+ \quad (47)$$

Since by (44) matrices

$$\begin{bmatrix} \mathbf{A} - \mathbf{LC} & \mathbf{F} \\ -\mathbf{C} & \mathbf{I}_p \end{bmatrix}, \begin{bmatrix} \mathbf{A} + \mathbf{KC} & \mathbf{0} \\ -\mathbf{C} & \mathbf{\Phi} \end{bmatrix} \quad (48)$$

of size $((n+p) \times (n+p))$ are similar, we have:

$$rank \begin{bmatrix} \mathbf{I}_n z - \mathbf{A} + \mathbf{LC} & -\mathbf{F} \\ \mathbf{C} & (z-1)\mathbf{I}_p \end{bmatrix} = rank \begin{bmatrix} \mathbf{I}_n z - \mathbf{A} - \mathbf{KC} & \mathbf{0} \\ \mathbf{C} & \mathbf{I}_p z - \mathbf{\Phi} \end{bmatrix} \quad (49)$$

From (47) and (49), we obtain:

$$rank \begin{bmatrix} \mathbf{I}_n z - \mathbf{A} + \mathbf{LC} & -\mathbf{F} \\ \mathbf{C} & (z-1)\mathbf{I}_p \end{bmatrix} = (n + p), \forall z \in C^+ \quad (50)$$

Furthermore, from (50) and Lemma 5 it follows that matrix

$$\begin{bmatrix} \mathbf{A} - \mathbf{LC} & \mathbf{F} \\ -\mathbf{C} & \mathbf{I}_p \end{bmatrix}$$

with **L** and **F** given by (42), (43) respectively is Schur stable and therefore according to (8) the system described by equations (4) and (5) with **L** and **F** given by (42) and (43), respectively is a proportional-integral observer of order *n* of the system described by equations (1) and (2). This completes the proof.

The sufficiency part of the proof of Theorem 1 provides a construction of the matrices **L** and **F** of proportional-integral observer of order *n* for the system described by equations (1) and (2). The major steps of this construction are given below.

**CONSTRUCTION**

*Given*: **A**, **B** and **C**

*Find*: **L** and **F**

*Step 1:* Check condition (a) of Theorem 1. If this condition is satisfied go to *Step 2*. If condition (a) of Theorem 1 is not satisfied, the construction of a proportional-integral observer of order *n* is impossible.

*Step 2:* Delectability of the pair (**A**, **C**) and Lemma 6 imply the existence of a matrix **K** over $\mathcal{R}$ of size *(n x p)* such that the matrix [**A**+**KC**] is Hurwitz stable. The matrix **K** can be calculated as in the proof of Lemma 6.

*Step 3:* Find a non-singular matrix **T** of size *(n x n)* such that

$$\mathbf{C} = [\mathbf{I}_p, \mathbf{0}]\mathbf{T}$$

*Step 4:* Let **Φ** be an arbitrary nonzero Schur stable matrix over $\mathcal{R}$ of size *(p x p)*. Further, let **Λ** be an arbitrary matrix over $\mathcal{R}$ of size *((n-p) x p)*. Put

$$\mathbf{X} = \mathbf{T}^{-1} \begin{bmatrix} (-\mathbf{\Phi} + \mathbf{I}_p) \\ \mathbf{\Lambda} \end{bmatrix}$$

$$\mathbf{L} = \mathbf{X} - \mathbf{K}$$

$$\mathbf{F} = -(\mathbf{A} - \mathbf{LC})\mathbf{X} + \mathbf{X}(-\mathbf{CX} + \mathbf{I}_p)$$

**5. CONCLUSIONS**

In this paper, by using basic concepts and basic results from linear systems and control theory as well as of the theory of matrices, the problem of the design of a full order proportional-integral observer for the state estimation of discrete-time linear time-invariant systems is studied and completely solved. The proof of the main results of this paper is constructive and furnishes a simple procedure for the construction of full order proportional-integral observer.